\documentclass[11pt]{article}
\usepackage{amsfonts}
\usepackage{amsthm,amsmath,amssymb}

\newtheorem{theorem}{Theorem}[section]
\newtheorem{conjecture}[theorem]{Conjecture}
\newtheorem{lemma}[theorem]{Lemma}

\theoremstyle{definition}

\newcommand{\reals}{{\mathbb R}}

\newcommand{\hilbert}{\bigcirc\kern -0.8em
              {\rm\scriptstyle {H}\;}}

\begin{document}
\renewcommand*\descriptionlabel[1]{{\bf\hspace\labelsep #1\hfil}}
\title{Slowdown for time inhomogeneous branching Brownian motion}
\author{Ming Fang\thanks{Mathematical Sciences Research Institute,
17 Gauss Way,  Berkeley, CA 94720-5070, fang0086@umn.edu} \and
Ofer Zeitouni\thanks{Weizmann institute
and University of Minnesota,
zeitouni@math.umn.edu.
Partially supported by
 NSF grant \#DMS-1106627 and a by
 grant from the Israel Science Foundation.
 }}
\maketitle
\abstract{We consider the maximal displacement of
one dimensional branching Brownian motion with
(macroscopically) time varying profiles. For monotone decreasing variances,
we show that the correction from linear displacement is not logarithmic but
rather proportional
to $T^{1/3}$. We conjecture that this is the worse case correction possible.}
\newcommand{\se}{ {\sigma}_{\mbox{\rm eff}}}
\section{Introduction and statement of results}
\setcounter{equation}{0}
\setcounter{figure}{0}
\setcounter{subsection}{1}

The classical branching Brownian motion (BBM)
model in $\reals$
can be described probabilistically as follows. At time $t=0$,
one particle exists and is located at the origin. This particle
starts performing Brownian motion, up to an exponentially distributed random
time. At that time, the particle instantaneously splits into two independent
particles, and those start afresh performing Brownian motion until their
(independent) exponential clock rings.

We introduce some notation. Let ${\cal N}_t$
denote the collection of particles
alive at time $t$, set $N(t)=|{\cal N}_t|$,
and for any particle $v\in {\cal N}_t$, let $x_v(s), s\in [0,t]$
denote the (Brownian) trajectory performed by the particle and its ancestors.
$N(t)$ is a continuous time branching process, and
it is straightforward to verify that $N(t)e^{-t}$ is a Martingale, which
converges almost surely to a positive, finite
random variable $n_\infty$. In particular, $t^{-1}\log N(t)$ converges
almost surely to $1$.

We will be interested in the location of the maximal particle, i.e.
in the random variable
$$ M_t=\max_{v\in {\cal N}_t} x_v(t)\,.$$
As is well known, the distribution function $F(x,t)=P(M_t \geq x)$
satisfies the Kolmogorov-Petrovskii-Piskunov equation (also
attributed to Fisher)
$$ \frac{\partial F}{\partial t}(x,t)=\frac12 \frac{\partial^2 F}{\partial^2 x}
(x,t)+F(x,t)(1-F(x,t))\,, \quad
F(x,0)={\bf 1}_{x\leq 0}\,.$$
See \cite{M75} for a probabilistic interpretation of the KPP equation.

In a seminal work, Bramson \cite{Br78} showed among other facts that
\begin{equation}
  \label{eq-of1}
  M_t=m(t)+O_P(1)\,,\quad
  m(t)=\sqrt{2} t-\frac{3}{2} \cdot \frac{1}{\sqrt{2}} \log t\,,
\end{equation}
in the sense that for any $\epsilon$ there is a $K_\epsilon$ so that
$$P(|M_t-m(t)|>K_\epsilon)\leq \epsilon\,.$$
In particular, $\textrm{Med}(M_t)=m(t)+O(1)$, where $\textrm{Med}(M_t)$ denotes
the median of $M_t$.
(In subsequent work \cite{Br83}, Bramson also discusses convergence to a shifted
traveling wave, but this is not the focus of the current work.) Analogues of
\eqref{eq-of1} also hold in the setup of discrete time branching random walks,
see \cite{ABR}; for a recent convergence result for BRWs, see \cite{Aid}.

%Ming
The leading term in \eqref{eq-of1}, linear in time,
%Ming
is a relatively
straight--forward consequence of large deviations computations and
the first and second moment methods; in particular, the coefficient
$\sqrt{2}$ would be the same if instead of BBM, one would consider the maximum
of $e^t$ independent Brownian motions run for time $t$. On the other hand, the
logarithmic correction term in \eqref{eq-of1}
is more subtle, and reflects the correlation structure of the BBM: for
the maxima of independent BMs, the $3/2$ multiplying the logarithmic term is
replaced by $1/2$. For a pedestrian introduction to these issues, see the
lecture notes \cite{Zei12}.

Our goal in this paper is to address situations in which the diffusivity of
the Brownian motion changes in time, in a macroscopic scale. This is
motivated in part by our earlier work \cite{FangZeitouni}, in which we showed
that the corrrection factor $3/2$ multiplying the logarithmic term
can be replaced by different, and eventually much larger, factors.
This naturally leads to the question, whether larger-than-logarithmic
corrections are possible. Our goal here is to answer this question in the
affirmative.

We turn to the description of the time inhomogeneous BBM that we consider.
Fix $T$ (eventually, large). We consider
the
BBM model where at time $t$, all particles
move independently as Brownian motions with variance
$\sigma_T^2(t)=\sigma^2(t/T)$, and
branch independently at rate $1$. Here, $\sigma$ is a smooth,
strictly decreasing
function
on $[0,1]$ with range in a compact subset of $(0,\infty)$, whose derivative is
bounded above by a strictly negative constant.
Define ${\cal N}_t$, $\{x_v(t)\}_{v\in {\cal N}_t}$ and $M_t$, $t\in [0,T]$,
as in the case of
time homogeneous BBM.
This model has been considered before in \cite{DS}.
%Let $\xi_i(t), i=1,\ldots,N_t$ denote the set
%of particles at time $t$, let
%$M_t=\max_{i} \xi_i(t)$ denote the maximal displacement of the BBM.
%We are interested in $M_T$.
Our main result is the following.
\begin{theorem}
  \label{theo-n13}
  With notation as above,
  %$\sigma(\cdot)$ as in \eqref{eq-sigma1},
  we have that
  \begin{equation}
    \label{eq-dispn13}
    \textrm{Med}(M_T)=v_\sigma T-g_{\sigma}(T)\,,
  \end{equation}
  where $v_\sigma$ is defined in \eqref{eq-vsn13}, and
  \begin{equation}
    \label{UB}
    0<\liminf_{T\to\infty}\frac{ g_\sigma(T)}{T^{1/3}}
    \leq \limsup_{T\to\infty}\frac{ g_\sigma(T)}{T^{1/3}}
    <\infty
  \end{equation}
\end{theorem}
We emphasize that it is already known a-priori \cite{Fang}
that
$\{M_T-\textrm{Med}(M_T)\}$ is
a tight sequence;
in fact, the tails estimates
in \cite{Fang} are
strong enough to allow one to replace, both in the statement above and
in Theorem \ref{theo-n13}, the median $\textrm{Med}(M_T)$ by $EM_T$.

\section{Proofs}
Before bringing the proof of Theorem \ref{theo-n13}, we collect some
preliminary information concerning the path of individual particles. With
$W_\cdot$ and $\tilde W_\cdot$ denoting standard Brownian
motions,
let
$$ X_t=\int_0^t \sigma_T(s) dW_s\,,\quad t\in [0,T]. $$
Let $\tau(t)=\int_0^t \sigma_T^2(s)ds$. Clearly, $X_\cdot$ has the same law
as $\tilde W_{\tau(t)}$. The following is a standard adaptation  of Schilder's
theorem \cite{Schilder,DZ}, using the scaling properties of Brownian motion.
\begin{theorem}[Schilder]
  \label{thm-schilder}
  Define $Z_t=  \frac{1}{T} X_{t/T}, t\in [0,1]$. Then $Z_t$ satisfies
a large deviation principle in $C_0[0,1]$ of speed $T$ and rate function
$$I(f)=\left\{\begin{array}{ll}
  \int_0^1 \frac{f'(s)^2}{2\sigma^2(s)} ds\,,& f\in H_1[0,1]\,,
  \\
  \infty\,,& \mbox{\rm else}
\end{array}
  \right.\,.$$
\end{theorem}
Here, $H_1[0,1]$ is the space of absolutely continuous function
on $[0,1]$ that vanish at $0$, whose (almost everywhere
defined)
derivative is
square-integrable.
%{\textit Proof of Theorem \ref{thm-schilder}}

We now wish to define a barrier for the particle systems that is
unlikely to be crossed. This barrier will also serve as a natural
candidate for a change of measure. Recall that at time $t$, with overwhelming
probability there are at most $e^{t+o(t)}$ particles alive in the
system. Thus, it becomes unlikely that any particle crosses
a boundary of the form $Tf(\cdot/T)$ if, at any time,
  $$J_t(f):=\int_0^t \frac{f'(s)^2}{2\sigma^2(s)} ds>t\,.$$
  This motivates the following lemma.
  \begin{lemma}
    \label{lem-varsol}
    Assume $\sigma$ is strictly decreasing. Then the solution of the
    variational problem
    \begin{equation}\label{eq-vsn13a}
     v_\sigma:= \sup \{f(1): J_t(f)\leq t, t\in [0,1]\}
   \end{equation}
exists, and the unique minimizing path is the function
\begin{equation}
  \label{eq-varsol}
  \bar f(t)=\sqrt{2}\int_0^t \sigma(s) ds\,.
\end{equation}
In particular,
    \begin{equation}\label{eq-vsn13}
      v_\sigma=\sqrt{2}\int_0^1 \sigma(s) ds\,.
    \end{equation}
\end{lemma}
\textit{Proof of Lemma \ref{lem-varsol}:}
We are going to prove that no other functions can do better than $\bar{f}$. That is, if some absolutely continuous function $g$ satisfies $g(0)=0$ and the constraint $J_t(g)\leq t$ for all $0\leq t\leq 1$, then $g(1)\leq \bar{f}(1)=v_{\sigma}$. In fact, denote $\phi(t)=J_t(g)\leq t$ for $0\leq t\leq 1$, and then $\phi'(t)=\frac{g'(t)^2}{2\sigma^2(t)}$ a.e.. We can write $g^2(1)$ as
\begin{eqnarray*}
  g^2(1)&=&\left(\int_0^1g'(t)dt\right)^2= \left(\int_0^1 \sqrt{2\phi'(t)}\sigma(t)dt\right)^2.
\end{eqnarray*}
Using H\"{o}lder's inequality, we have
\begin{eqnarray*}
 g^2(1) \leq 2\left( \int_0^1 \phi'(t)\sigma(t)dt\right) \left( \int_0^1\sigma(t)dt \right)
  = \sqrt{2}v_{\sigma}\left( \int_0^1 \phi'(t)\sigma(t)dt\right).
\end{eqnarray*}
Using integration by parts, the above is equal to
\begin{eqnarray*}
  \sqrt{2}v_{\sigma}\left( \phi(1)\sigma(1)-\int_0^1 \phi(t)\sigma'(t)dt\right).
\end{eqnarray*}
Since $\phi(t)\leq t$ and $\sigma'(t)\leq 0$ for all $0\leq t\leq 1$, the above is less than or equal to
\begin{eqnarray*}
  \sqrt{2}v_{\sigma} \left( \sigma(1)-\int_0^1 t\sigma'(t)dt\right)= \sqrt{2}v_{\sigma}\int_0^1\sigma(t)dt=v_{\sigma}^2,
\end{eqnarray*}
where we apply integration by parts in the first equality. This completes the proof.
\qed

%By strict convexity of $J_t(\cdot)$, a solution $g$ of the variational
%problem exists and is unique, and
%it satisfies $g'(s)\geq 0$ for all $s\in [0,1]$. Note that $\bar f$ is
%the solution of the differential equation
%$\bar f'(s)=\sqrt{2}\sigma(s)$. Fix $\delta>0$ so that with
%$\theta_1=\inf\{t>0: J_t( g)<t-\delta\}\wedge  1$ and
%$\theta_2=\inf\{t>\theta_1: J_t( g)= t-\delta/2\}\wedge 1$,
%one has $\theta_1<1$.
%(Such a $\delta$ exists if $\bar f\neq g$.)
%Let $A=\{s\in [\theta_1,\theta_2]:  g'(s)<\sqrt{2}\sigma(s)\}$
%and  $B=\{s\in [\theta_1,\theta_2]:  g'(s)>\sqrt{2}\sigma(s)\}$.
%If $\bar f\neq g$ then
%%reducing $\delta$ if necessary,
%the Lebesgue measures of $A, B$ are positive, and furthermore
%there is an $\epsilon_0>0$ so that for each $\epsilon\in (0,\epsilon_0)$ there
%exist sets $A_1=A_1(\epsilon)\subset A$
%and $B_1=B_1(\epsilon)\subset B$, both of equal
%positive Lebesgue measure, so that
%$x<y$, $g'(x)<(\sqrt{2}-\epsilon)\sigma(x)$,
%$g'(y)>(\sqrt{2}+\epsilon)\sigma(y)$
%for all $x\in A_1,y\in B_1$. The function $\bar g$ obtained
%from $g$ by setting $\bar g(0)=0$ and
%$$\bar g'(x)=\left\{\begin{array}{ll}
%  g'(x)+\epsilon \sigma(x) \,, &x\in A_1\,,\\
%  g'(x)-\epsilon \sigma(x)\,, & x\in B_1\,,\\
%  g'(x)\,, & \mbox{\rm otherwise}\,,
%\end{array}
%\right.$$
%satisfies, for $\epsilon$ small enough,
%the constraint  $J_t(\bar g)\leq t$
%(here we used that $A_1,B_1\in (\theta_1,\theta_2)$.)
%Further, due to the monotonicity of
%$\sigma(\cdot)$,
%$\bar g(1)\geq g(1)$.
%This contradicts the uniqueness of the minimizer.
%Hence, necessarily, $g=\bar f$, as claimed.
%\qed
%

\noindent
\textit{Proof of Theorem \ref{theo-n13}:}
We begin with the upper bound in
%It remains  to provide a proof of the upper bound in
\eqref{eq-dispn13}. The first step is to show that
in fact, no particle will be found significantly above
$T\bar f(t/T)$.
\begin{lemma}
  \label{lem-ubn13}
  There exists $C$ large enough such that, with
  $${\cal A}=\{
  \exists t\in [0,T], v\in {\cal N}_t,
  %i\in \{1,\ldots, N_t\}
    x_v(t)>T\bar f(t/T)+C\log T\}\,,$$
    it holds that
  \begin{equation}
    \label{eq-lemubn13}
    P({\cal A})
    \to_{T\to\infty}
    0\,.
  \end{equation}
\end{lemma}
\textit{Proof of Lemma \ref{lem-ubn13}:}
Recall the process $X_\cdot$ in $C_0[0,T]$, whose law we denote by $P_0$.
Consider the change of measure
with Radon--Nykodim derivative
%Ming
\begin{eqnarray}
  \label{eq-RNnew}
  \frac{dP_1}{dP_0}|_{ {\cal F}_t}&=&\exp\left(-\int_0^t\frac{\bar
f'(s/T)}{\sigma^2(s/T)} dX_s
-\frac12\int_0^t \frac{(\bar
f'(s/T))^2}{\sigma^2(s/T)} ds\right)\nonumber\\
&=&\exp\left(-\int_0^t\frac{\sqrt2}{\sigma(s/T)} dX_s
-t\right)\,.
\end{eqnarray}
The process $X_\cdot$ under $P_0$ is the same
as
%the law of
the process $X_\cdot+T\bar f(\cdot/T)$ under $P_1$.
Note that for any $t\leq T$,
%Ming
\begin{equation}
  \label{eq-after13}
  \int_0^t \frac{\sqrt2}{\sigma(s/T)} dX_s=
  \frac{\sqrt{2}X_t}{\sigma(t/T)}+\frac{\sqrt{2}}{T}
  \int_0^t X_s \frac{\sigma'(s/T)}
  {\sigma^2(s/T)}ds\,.
\end{equation}
We then have,
with $\tau=\inf\{t\leq T: X_t\geq C\log T\}$,
on the event $\tau\leq T$,
%Ming
\begin{eqnarray*}
\int_0^{\tau} \frac{(\bar f'(s/T))}{\sigma^2(s/T)} dX_s
&\geq &\frac{\sqrt{2}C\log T}{\sigma(t/T)}+\frac{\sqrt{2}C\log T}{T}
  \int_0^t  \frac{\sigma'(s/T)}
  {\sigma^2(s/T)}ds\\
  &= &\frac{\sqrt{2}C\log T}{\sigma(0)}\,,
\end{eqnarray*}
%Ming
and therefore, with $\tau'=\inf\{t\leq T: X_t\geq T\bar f(t/T)+C\log T\}$, we have, for $k\leq T$,
\begin{eqnarray*}
%&&P_0(\exists t\leq T:
%X_t\geq T\bar f(t/T)+C\log T)=
%P_1(\exists t\leq T:
%X_t\geq C\log T)\\
P_0(\tau'\in[k-1,k))&=& P_1(\tau\in[k-1,k))=
E_{P_0}\left(\frac{dP_1}{dP_0} {\bf 1}_{\tau\in[k-1,k)}\right)\\
&\leq &
E_{P_0}\left({\bf 1}_{\tau\in[k-1,k)}
\exp\left(-\frac{\sqrt{2}C\log T}{\sigma(0)}-\tau\right)\right).
%E_{P_0}\left(\exp\left(-C \log T/\sqrt{2}+C_1-
\end{eqnarray*}

Define
%Thus,
%using the first moment method, letting
$$\theta=\inf\{t\leq T: \mbox{\rm there is $v\in {\cal N}_T$
so that $x_v(t)\geq
T\bar f(t/T)+C\log T$}\}\,,$$
%Ming
and $Z_k$ to be the number of particles $z\in{\cal N}_k$ such that $x_v(t)\leq T\bar f(t/T)+C\log T$ for all $t\leq k-1$ and $x_v(t)\geq T\bar f(t/T)+C\log T$ for some $k-1\leq t\leq k$.
Then,
$$P(\theta\leq T)\leq \sum_{k=1}^T
P(\theta\in [k-1,k))\leq P(Z_k\geq 1)\,,$$
and, using a first moment computation, we obtain
$$P(Z_k\geq 1)\leq EZ_k \leq e^k P_0(\tau'\in[k-1,k)) \leq \exp\left(-\frac{\sqrt{2}C\log T}{\sigma(0)}+1\right).$$
Therefore,
$$P(\theta\leq T)
\leq
T\exp\left(-\frac{\sqrt{2}C\log T}{\sigma(0)}+1\right).$$
This completes the proof of Lemma
  \ref{lem-ubn13}.
  \qed

We need one more technical estimate.
\begin{lemma}
  \label{lem-almostfinaln13}
%  Fix $\epsilon>0$.
With $X_\cdot$ and $C$ as in Lemma
\ref{lem-ubn13}, there exists a constant $C'\in (0,1)$ so that
\begin{equation}
  \label{eq-finaln13}
  e^T
  P_0(X_t\leq T\bar f(t/T)+C\log T, t\in[0,T], X_T\geq T\bar f(1)-C'T^{1/3})
  \to_{T\to\infty} 0\,.
\end{equation}
\end{lemma}
\textit{Proof of Lemma \ref{lem-almostfinaln13}:}
Fix $C'\in (0,1)$.
We apply a change of measure similar to the one used
in Lemma \ref{lem-almostfinaln13}, whose notation we continue to use.
We deduce
the existence of  positive constants $c_1,c_2$ (independent
of $T$) such that
%Ming
\begin{eqnarray*}
  &&P_0(X_t\leq T\bar f(t/T)+C\log T, t\in[0,T], X_T\geq T \bar f(1)-C'T^{1/3})
\\&\leq &
e^{-T}e^{c_1(C'T^{1/3}+\log T)}\\
&&\cdot
E_{P_0}\left(\exp\left(\frac{c_2}{T} \int_0^T X_s
ds\right){\bf 1}_{X_T\geq - C'T^{1/3}}
{\bf 1}_{X_t\leq 0, t\leq T}\right),
\end{eqnarray*}
where here we used that $-\sigma'$ is bounded below by a positive constant and
$\sigma$ is bounded above.
%for the first time the explicit affine form of $\sigma$.
By representing $X_\cdot$ as a time-changed Brownian motion, the lemma will
follows (for a small enough $C'$)
if we can show that for any constant $c_3$ there
exists a
%Ming
$c_4=c_4(c_3)>0$ independent of $C'\in (0,1)$ such that
\begin{equation}
  \label{eq-bmn13}
  {\tt D}:=  E\left(\exp\left(\frac{c_3}{T} \int_0^T B_s
ds\right){\bf 1}_{B_T\geq - C'T^{1/3}}
{\bf 1}_{B_t\leq 0, t\leq T}\right)\leq e^{-c_4T^{1/3}}\,,
\end{equation}
where $\{B_t\}_{t\geq 0}$ is a Brownian motion
started at $-C\log T$.
Note however that
$${\tt D}\leq
   E\left(\exp\left(-\frac{c_3}{T} \int_0^T |B_s|
ds\right){\bf 1}_{|B_T|\leq  T^{1/3}}
\right)e^{c_5 \log T}\leq e^{-c_4T^{1/3}}\,,
$$
where here $B_\cdot$ is a Brownian motion
started at $0$ and
the last inequality is a consequence of known estimates
for Brownian motion, see
e.g. \cite[Formula 1.1.8.7, pg. 141]{BS}.
\qed

%Ming
{\bf Remark}\quad
The estimate in
  \eqref{eq-bmn13} can also be derived probabilistically. Here
  is a sketch.
  It is clearly enough to estimate the expectation on the event
  ${\tt F}:=\left\{\mbox{\rm Leb}(\{s: 0\geq B_s\geq -c_5 T^{1/3}\}) \geq 1/2\right\}$.
  But $P({\tt F})$ decays exponentially in $T^{1/3}$. Some more details
  are provided in
  \cite{Zei12}.

We have completed all steps required for the proof of the upper
bound in Theorem
\ref{theo-n13}.
%Indeed, it only remains to prove the upper bound in the Theorem.
Due to the strong tightness result in \cite{Fang} and Lemma
\ref{lem-ubn13}, it
is enough to show that
$$P(\{M_T\geq \bar Tf(1)-C'T^{1/3}\}\cap {\cal A}^\complement)\to 0\,.$$
This follows from the first moment method and Lemma
  \ref{lem-almostfinaln13}.

  We turn to the proof of the lower bound.
%Consider the path of particles at time $T$, denoted by
%$x_i(t)$, where $i=1,\ldots,N(T)$ and $0\leq t\leq T$.
Call a particle $v\in {\cal N}_T$
\textit{good}
if
$$T\bar{f}(t/T)-T^{1/3}\leq x_v(t)\leq T\bar{f}(t/T), \text{ for all } t\leq T.$$
Set
$${\cal M}=\sum_{v\in {\cal N}_T}
{\bf 1}_{\textrm{$v$ is a good particle}}\,.$$

\begin{lemma}
  \label{lem-ming1}
  There exists a constant $C>0$ such that
  $$P({\cal M}\geq 1)\geq e^{-CT^{1/3}}.$$
\end{lemma}
\textit{Proof of Lemma \ref{lem-ming1}:}
Recall the process $X_\cdot$ in $C_0[0,T]$, whose law we denoted by $P_0$, and
the measure $P_1$ defined by \eqref{eq-RNnew}.
%Consider the change of measure
%with Radon--Nykodim derivative
%\begin{proof}
%  We use a change of measure and second moment method. Consider the change of measure
%with Radon--Nykodim derivative
%\begin{eqnarray*}
%\frac{dP_1}{dP_0}&=&\exp\left(\int_0^T\frac{\bar
%f'(s/T)}{\sigma^2(s/T)} dX_s
%-\frac12\int_0^T \frac{(\bar
%f'(s/T))^2}{\sigma^2(s/T)} ds\right)\\
%&=& \exp\left(\int_0^T\frac{\bar
%f'(s/T)}{\sigma^2(s/T)} d\tilde{X}_s
%+\frac12\int_0^T \frac{(\bar
%f'(s/T))^2}{\sigma^2(s/T)} ds\right)\\
%&=& \exp\left(\int_0^T\frac{1}{\sigma(s/T)} d\tilde{X}_s
%+T\right)\,.
%\end{eqnarray*}
%Here $X_.$ under $P_0$ has the same law as $\tilde{X}_.=X_.-T\bar{f}(\cdot/T)$ under $P_1$.
%%
We then calculate the first moment
\begin{eqnarray*}
  E{\cal M}&=&E\sum_{v\in {\cal N}(T)}
  {\bf 1}_{\textrm{$v$ is a good particle}}\\
&=& e^TP_0(T\bar{f}(t/T)-T^{1/3}\leq X_t\leq T\bar{f}(t/T), \text{ for all } t\leq T)\\
&=&E_{P_0}\left[\exp\left(-\int_0^T\frac{\sqrt{2}}{\sigma(s/T)} d{X}_s\right)
{\bf 1}_{\{-T^{1/3}\leq {X}_t\leq 0, \text{ for all } t\leq T\}}\right].
\end{eqnarray*}
%Ming
Repeating the computation in \eqref{eq-after13}, we conclude that
%By integration by parts, one can rewrite the integral in the above expression as
%$$-\int_0^T\frac{1}{\sigma(s/T)} d\tilde{X}_s = -\frac{\tilde{X}_T}{\sigma(1)} - \frac{1}{T}\int_0^T\frac{\tilde{X}_s\sigma'(s/T)}{\sigma^2(s/T)}ds.$$
%Recall that $\sigma'(\cdot)$ is negative. In the event of $\{-T^{1/3}\leq \tilde{X}_t\leq 0, \text{ for all } t\leq T\}$, the above is greater than or equal to
%$$T^{1/3} \left(\frac{1}{T}\int_0^T\frac{\sigma'(s/T)}{\sigma^2(s/T)}ds\right) \geq T^{1/3}\min_{0\leq x\leq 1}\frac{\sigma'(x)}{\sigma^2(x)}:=-c_1T^{1/3}$$
%where $c_1=-\min_{0\leq x\leq 1}\frac{\sigma'(x)}{\sigma^2(x)}>0$.
%The above argument gives a lower bound for the first moment $E{\cal N}$ as follows
$$E{\cal M} \geq \exp\left(-c_6T^{1/3}\right) P_0\left(-T^{1/3}\leq {X}_t\leq 0, \text{ for all } t\leq T\right).$$
Since under $P_0$, $X_\cdot$ is a time changed Brownian motion, we have that
%Since $\tilde{X}_t$ under $P_1$ is a time changed Brownian motion, a M\"{o}gul'skii result shows that
$$P_0\left(-T^{1/3}\leq {X}_t\leq 0, \text{ for all } t\leq T\right)
\geq e^{-c_7T^{1/3}}$$
for some $c_7>0$.
Hence,
$$E{\cal M}\geq e^{-c_8T^{1/3}}$$
for some $c_8>0$.

We next
derive an upper bound for the second moment $E{\cal M}^2$. By definition,
\begin{eqnarray*}
  E{\cal M}^2&=&E\sum_{v,v'\in {\cal N}_T} {\bf 1}_{\textrm{$v,v'$ are
  good particle}}.
\end{eqnarray*}
When $v\neq v'$, we
let $t_{vv'}$ be the branching time of the last
common ancestor of $v$ and $v'$.
Then, the paths $\{x_v(s)\}_{0\leq s\leq T}$ and
$\{x_{v'}(s)-x_{v'}(t_{vv'})\}_{t_{vv'}\leq s\leq T}$
are independent. Applying a change of measure similar to that
used in the computation of $E{\cal M}$,
we can bound above
\begin{eqnarray*}
 && E{\cal M}^2 \leq E{\cal M}\\
  &&+\int_0^T e^{2T-t}
  E\left[\exp\left(-\int_0^T\frac{\sqrt{2}}{\sigma(s/T)} dX^1_s-T\right)
  {\bf 1}_{\{-T^{1/3}\leq X^1_s\leq 0, \text{ for all } 0\leq s\leq T\}}\right.\\
  && \;\;\;\cdot\left.\exp\left(-\int_t^T\frac{\sqrt{2}}{\sigma(s/T)} dX^2_s-(T-t)\right)
  {\bf 1}_{\{-T^{1/3}\leq X^2_s\leq T^{1/3}, \text{ for all } t\leq s\leq T\}}\right]dt,
\end{eqnarray*}
where $X^1_.$ and $X^2_.$ are two i.i.d. copies of $X_.$
(under the law $P_0$). The above is equal to
\begin{eqnarray*}
  && E{\cal M}+\int_0^T
  E\left[\exp\left(-\int_0^T\frac{\sqrt{2}}{\sigma(s/T)} dX^1_s\right)
  {\bf 1}_{\{-T^{1/3}\leq X^1_s\leq 0, \text{ for all } 0\leq s\leq T\}}\right.\\
  && \;\;\;\cdot\left.\exp\left(-\int_t^T\frac{\sqrt{2}}{\sigma(s/T)} dX^2_s\right)
  {\bf 1}_{\{-T^{1/3}\leq X^2_s\leq T^{1/3}, \text{ for all } t\leq s\leq T\}}\right]dt.
\end{eqnarray*}

%We give the desired upper bound for the second moment $E{\cal N}^2$.
On the event $\{-T^{1/3}\leq X^1_s\leq 0, \text{ for all }
0\leq s\leq T\}$, using integration by parts, one has
$$-\int_0^T\frac{1}{\sigma(s/T)} dX^1_s = -\frac{X^1_T}{\sigma(1)} - \frac{1}{T}\int_0^T\frac{X^1_s\sigma'(s/T)}{\sigma^2(s/T)}ds\leq \frac{T^{1/3}}{\sigma(1)}.$$
Similarly, on the event $\{-T^{1/3}\leq X^2_s\leq T^{1/3},
\text{ for all } t\leq s\leq T\}$, using integration by parts, one has
%Ming
\begin{eqnarray*}
-\int_t^T\frac{1}{\sigma(s/T)} dX^2_s
& =& -\frac{X^2_T}{\sigma(1)}+\frac{X^2_t}{\sigma(t/T)} - \frac{1}{T}\int_t^T\frac{X^2_s\sigma'(s/T)}{\sigma^2(s/T)}ds\\
&\leq & \frac{3T^{1/3}}{\sigma(1)}\,.
\end{eqnarray*}
%Ming
Therefore, $E{\cal M}^2$ is bounded above by $E{\cal M}$ plus
%Ming
\begin{eqnarray*}
&&\int_0^T e^{4T^{1/3}/\sigma(1)}\\
&&
  E\left[
  {\bf 1}_{\{-T^{1/3}\leq X^1_s\leq 0, \text{ for all } 0\leq s\leq T\}}
  {\bf 1}_{\{-T^{1/3}\leq X^2_s\leq T^{1/3}, \text{ for all } t\leq s\leq T\}}\right]dt\,.
\end{eqnarray*}
The latter integral is
less than or equal to
%Ming
$$\int_0^T e^{4T^{1/3}/\sigma(1)}dt\leq e^{c_9T^{1/3}}$$
for some $c_9>0$.
Hence,
one can apply the second moment method:
$$P({\cal M}\geq 1)=P({\cal M}>0)\geq \frac{(E{\cal M})^2}{E{\cal M}^2}\geq
e^{-c_{10}T^{1/3}}.$$
This completes the proof of the lemma by letting $C=c_{10}>0$.
\qed

By a direct first moment computation (or using the already proved upper
bound \eqref{UB}),
the minimum position of particles at time $AT^{1/3}$ is, with high probability,
greater than or equal to $-A'T^{1/3}$ for a constant $A'$ depending on $A$.
Choosing $A>C$ (with $C$ as in Lemma \ref{lem-ming1}), and using the
independence of the motion of descendents of different particles in
${\cal N}_{AT^{1/3}}$,
the lower bound in \eqref{UB} follows.
  \qed

  \section{Discussion}
  Our choice of considering strictly decreasing diffusivity is not accidental:
  the computations in \cite{FangZeitouni} hint that this should correspond to a
  worse-case situation. In fact, we conjecture the following.
  \begin{conjecture}
    Let $\sigma$ be a smooth function from $[0,1]$ to a compact
    subset of $(0,\infty)$. Then, \eqref{eq-dispn13}
    %and the upper bound in
    %\eqref{UB} continue to
    holds with $v_\sigma$ determined by \eqref{eq-vsn13a}
    and $|g_\sigma(T)|=O(T^{1/3})$.
  \end{conjecture}

  In a slightly more technical direction, it would of course be of
  interest, in the setup of Theorem \ref{theo-n13},
  to show that $g_\sigma(T)/T^{1/3}$ converges as $T\to\infty$,
  and to evaluate the limit. Our methods are not refined enough to
  allow for that.

  Finally, we mention that results for the homogenization of the KPP
  equation are available, see e.g.  \cite{nolen,nolen1} and references therein.
  In the terminology we employ here, those results correspond to
  fast varying, or \textit{microscopic}, time inhomogenuities.


\begin{thebibliography}{99}
    \bibitem[ABR09]{ABR}
L. Addario-Berry and B.   Reed,  Minima in branching random walks,
 {\it Annals Probab.}  {\bf 37}  (2009),  pp.  1044--1079.
\bibitem[Ai11]{Aid} E. Aidekon, Convergence in law of the minimum of a branching
  random walk. arXiv:1101.1810v3. To appear, {\it Annals Probab.}
 \bibitem[BS96]{BS} A. N. Borodin and O. Salminen,
   {\it Handbook of Brownian motion - facts and formulae},
   Birkhauser, Basel (1996).
 \bibitem[Br78]{Br78} M.  Bramson,  Maximal displacement
of branching Brownian motion, {\it Comm. Pure Appl. Math.}
{\bf 31} (1978),
pp. 531--581.
\bibitem[Br83]{Br83} M. Bramson, Convergence of solutions of the Kolmogorov
equation to travelling waves, {\it Mem. Amer. Math. Soc.} {\bf 44} (1983),
No. 285.

\bibitem[BZ09]{BZ} M. Bramson and O. Zeitouni,
  Tightness for a family of recursion equations,
  {\it Annals Probab.} {\bf 37} (2009), pp. 615--653.
\bibitem[DZ98]{DZ} A. Dembo and O. Zeitouni,
{\it Large deviation techniques and applications}, 2nd edition,
Springer, New-York (1998).
\bibitem[DS88]{DS}
B. Derrida, B. and H.
Spohn,   Polymers on disordered trees, spin glasses, and traveling waves.
%New directions in statistical mechanics (Santa Barbara, CA, 1987).
{\it J. Statist. Phys.} {\bf   51}  (1988),  pp.  817--840.
\bibitem[Fa10]{Fang} M. Fang,
  Tightness for maxima of generalized branching random walks (2010).
  	arXiv:1012.0826v1 (2010).
      %\bibitem[FZ10]{FZ10} M. Fang and O. Zeitouni,
%		Consistent Minimal Displacement of Branching Random Walks,
%		{\it ECP} {\bf  15} (2010), pp. 106--118.
\bibitem[FZ11]{FangZeitouni} M. Fang and O. Zeitouni,
Branching Random Walks in Time Inhomogeneous Environments,
	arXiv:1112.1113v1 (2011).
\bibitem[M75]{M75} H. McKean, Application of Brownian motion to the equation of
Kolmogorov-Petrovskii-Piskunov, {\it Comm. Pure Appl. Math.} {\bf 28}
(1975), pp. 323-331.
\bibitem[N11]{nolen} J. Nolen, An invariance principle for
  random traveling waves in one dimension,
  {\it SIAM J. Math. Anal.} {\bf 43} (2011), pp. 153--188.
  %{\it Networks and Heterogeneous Media} {\bf 6} (2011), pp. 167--194.
\bibitem[NRRZ12]{nolen1} J. Nolen, J.-M. Roquejoffre, L. Ryzhik and A. Zlatos,
 Existence and non-existence of Fisher-KPP transition fronts,
 {\it Arch. Rat. Mech. Anal.} {\bf 203} (2012), pp. 217--246.
\bibitem[Ro11]{roberts} M. Roberts,
  A simple path to asymptotics for the frontier of a branching Brownian motion,
  to appear, {\it Annals Probab.} (2012). arXiv:1106.4771v4.
\bibitem[S66]{Schilder} M. Schilder, Some asymptotic formulae for Wiener integrals, {\it Trans. Amer. Math. Soc.}
  {\bf 125} (1966), pp. 63--85.
\bibitem[Z12]{Zei12} O. Zeitouni, Lecture notes on Branching random walks and the
  Gaussian free field, in preparation (2012).
\end{thebibliography}
  \end{document}